\documentclass[11pt]{article}

\usepackage[latin1]{inputenc} 
\usepackage{amsthm,amsmath,amssymb} 
\usepackage{graphicx}
\usepackage{color}
\usepackage{verbatim}

\usepackage{enumitem} 
\setlist[enumerate]{topsep=0pt,itemsep=-1ex,partopsep=1ex,parsep=1ex}

\theoremstyle{plain}
\newtheorem{theo}{Theorem}[section]

\newtheorem{lemma}[theo]{Lemma}

\theoremstyle{definition}

\newtheorem{alg}[theo]{Algorithm}
\newtheorem{claim}[theo]{Claim}

\newcommand{\mb}[1]{\mathbb{#1}}
\newcommand{\nib}[1]{\noindent {\bf #1}}

\newcommand{\sub}{\subseteq}

\newcommand{\sm}{\setminus}

\newcommand{\eps}{\varepsilon}
\newcommand{\es}{\emptyset}

\newcommand{\dD}{\delta}

\newcommand{\tT}{\theta}

\title{On the number of symbols that forces a transversal}

\author{Peter Keevash\thanks{Mathematical Institute, University of Oxford, Oxford, UK. 
E-mail: keevash@maths.ox.ac.uk.\newline \hspace*{1.5em}
Research supported in part by ERC Consolidator Grant 647678.}
\and Liana Yepremyan\thanks{Mathematical Institute, University of Oxford, Oxford, UK. 
E-mail: yepremyan@maths.ox.ac.uk.}}

\begin{document}
\maketitle

\begin{abstract}
Akbari and Alipour \cite{AA} conjectured that any Latin array of order $n$
with at least $n^2/2$ symbols contains a transversal.
We confirm this conjecture for large $n$, and moreover,
we show that $n^{399/200}$ symbols suffice.
\end{abstract}

\section{Introduction}

A Latin square of order $n$ is an $n$ by $n$ square 
with cells filled using $n$ symbols so that every
symbol appears once in each row and once in each column.
A transversal in a Latin square is a set of cells using
every row, column and symbol exactly once.
The study of transversals in Latin squares goes back to Euler in 1776;
his famous `36 officers problem' is equivalent to showing
that there is no Latin square of order 6 
that can be decomposed into 6 transversals
(this was finally solved by Tarry in 1900).
An even more fundamental question is whether
a Latin square always has a transversal.
A quick answer is `no', as shown by the addition table of $\mb{Z}_{2k}$,
but it remains open whether there is always a transversal 
when $n$ is odd (a conjecture of Ryser \cite{ryser}) 
or whether there is always a partial transversal size $n-1$
(a conjecture of Brualdi \cite{brualdi} and Stein \cite{stein}).
The best known positive result, due to Hatami and Shor \cite{HS},
is that there is always a partial transversal size $n-O(\log^2 n)$.

Given the apparent difficulty of finding transversals in Latin squares,
it is natural to ask if the problem becomes easier in Latin arrays with 
more symbols (now we fill a square with any number of symbols such
that every symbol appears at most once in each row and at most once in each column).
This problem was considered by Akbari and Alipour \cite{AA},
who conjectured that any Latin array of order $n$
with at least $n^2/2$ symbols contains a transversal.
Progress towards this conjecture was independently obtained by
Best, Hendrey, Wanless, Wilson and Wood \cite{BHWWW}
(who showed that $(2-\sqrt{2})n^2$ symbols suffice)
and Bar\'at and Nagy \cite{BN} 
(who showed that $3n^2/4$ symbols suffice).

We will henceforth adopt the standard graph theory translation of the problem,
where we consider a Latin array of order $n$ as a properly edge-coloured
complete bipartite graph $K_{n,n}$, with one part corresponding to rows, 
the other part to columns, and the colour of an edge is the symbol
in the corresponding cell of the array. In this language,
the Akbari-Alipour conjecture is that if there are at least $n^2/2$ colours
then there is a rainbow perfect matching.
Our main result confirms this conjecture for large $n$ in a strong form.

\begin{theo}\label{main}
Suppose the complete bipartite graph $K_{n,n}$ is properly edge-coloured
using $dn^2$ colours, where $n$ is sufficiently large and $d > n^{-1/200}$.
Then there is a rainbow perfect matching.
\end{theo}

Here the constant `200' could be somewhat improved,
but we have sacrificed some optimisations 
for the sake of readability of our proof,
as in any case the best bound we can obtain seems far from optimal
(it might even be true that $n^{1+o(1)}$ colours suffice!)

\section{Proof}

Here we give the proof of Theorem \ref{main},
assuming two lemmas that will be proved later in the paper.
Consider the complete bipartite graph $K_{n,n}$
with parts $A$ and $B$ both of size $n$,
and a proper edge-colouring using at least $dn^2$ colours,
where $n$ is sufficiently large and $d > n^{-1/200}$.

Let $G$ be any subgraph of $K_{n,n}$ obtained
by including exactly one edge of each colour.
Then $e(G) \ge dn^2$. We apply the following lemma,
which will be proved in the next section, to find a pair 
$(A_1,B_1)$ that satisfies Hall's condition for 
a perfect matching `robustly', so that it will still
satisfy Hall's condition after deleting small sets
of vertices from each part. Note that as $d > n^{-1/200}$
we obtain $|A_1|=|B_1| > n^{0.7}/3$. 

\begin{lemma} \label{robmatch}
There is $G_1=G[A_1,B_1]$ for some $A_1 \sub A$ and $B_1 \sub B$
with $|A_1|=|B_1| \ge d^{60} n/3$ such that 
$G_1$ has minimum degree at least $10^{-3} d|A_1|$,
and for any $S \sub A_1$ or $S \sub B_1$
we have $|N_{G_1}(S)| \ge \min\{2|S|,2|A_1|/3\}$.
\end{lemma}

We define a random subgraph $G^r$ of $K_{n,n}$ of `reserved colours' as follows.
Choose each colour independently with probability $p := n^{-0.32}$.
Let $G^r$ consist of all edges of all chosen colours. By Chernoff bounds, whp 
$|N_{G^r}(b) \cap A_1| = p|A_1| \pm (p|A_1|)^{2/3}$ and
$|N_{G^r}(b) \sm A_1| = p|A \sm A_1| \pm (p|A \sm A_1|)^{2/3}$ for all $b \in B$,
and similarly with $A$ and $B$ interchanged.
Let $G^* := (K_{n,n} \sm G^r) \sm (A_1 \cup B_1)$.
Then the minimum degree in $G^*$ satisfies
$\dD(G^*) \ge (1-p)(n-|A_1|) - (pn)^{2/3}$.

Let $M_2$ be a maximum size rainbow matching
in $G^* := (K_{n,n} \sm G^r) \sm (A_1 \cup B_1)$.
Let $A_2 = V(M_2) \cap A$ and $B_2 = V(M_2) \cap B$.
By a result of Gy\'arf\'as and S\'ark\"ozy \cite[Theorem 2]{GS}
we have $|A_2|=|B_2| \ge \dD(G^*) - 2 \dD(G^*)^{2/3} \ge (1-2p)(n-|A_1|)$,
as $pn = n^{0.68} \gg n^{2/3}$.

Let $G'_1$ be obtained from $G_1$ by deleting any edges
that use a colour used by $M_2$ and restricting
to some subsets $A'_1 \sub A_1$ and $B'_1 \sub B$
with $|A'_1|=|B'_1| = (1-10^{-4}d)|A_1|$ so that 
$G'_1$ has minimum degree at least $(10^{-3}-2\cdot 10^{-4})d|A_1|$.
To see that this is possible, note that we delete 
at most $n$ edges from $G'_1$, so each of $A_1$ and $B_1$
has at most $10^4 n/d|A_1| < 10^{-4}d|A_1|$ vertices at which 
we delete more than $10^{-4}d|A_1|$ edges.

\begin{figure}
\begin{center}
\includegraphics[scale = 0.6]{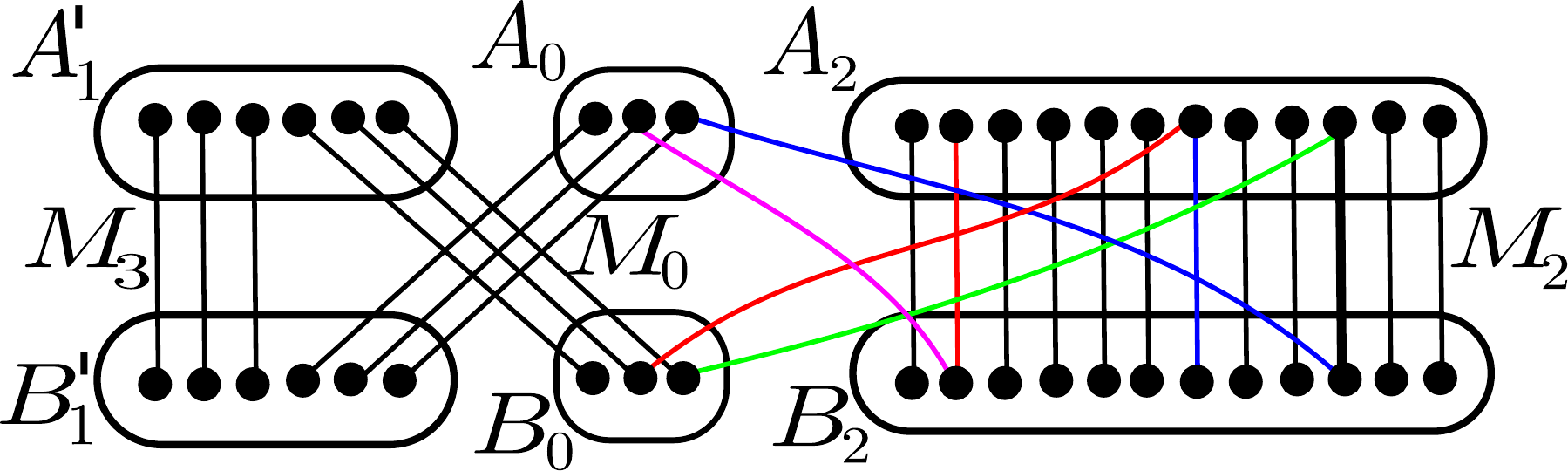}
\end{center}
\caption{Proof by picture}
\label{fig:pf}
\end{figure}

Let $A_0 = A \sm (A'_1 \cup A_2)$
and $B_0 = B \sm (B'_1 \cup B_2)$.
Note that $|A_0|=|B_0| < 2 \cdot 10^{-4} d|A_1|$,
as $pn = n^{0.68} \ll d|A_1|$.
The form of our intended rainbow matching
is illustrated by the black and/or vertical edges 
in Figure \ref{fig:pf} 
(the coloured diagonal edges illustrate
the augmentation algorithm used in Section \ref{sec:alg}).

Let $A'_0$ be the set of vertices $a$ in $A_0$
such that at least $|B'_1|/2$ of the edges
between $a$ and $B'_1$ have a colour used by $M_2$.
Define $B'_0$ similarly, interchanging $A$ and $B$.
We prove the following lemma in section \ref{sec:alg}.

\begin{lemma} \label{alg}
Both $A'_0$ and $B'_0$ have size at most $p|A_1|/4$.
\end{lemma}

Now we apply a greedy algorithm to construct 
a rainbow matching $M_2 \cup M_0$ where each
edge of $M_0$ joins $A_0 \cup B_0$ to $A'_1 \cup B'_1$.
We start by choosing these edges for vertices 
in $A'_0 \cup B'_0$ using colours in $G^r$. 
This is possible as the number of choices at each step 
is at least $p|A_1| - (p|A_1|)^{2/3} > 3p|A_1|/4$,
and at most $3 \cdot p|A_1|/4$ choices are forbidden 
due to using a colour or a vertex used at a previous step.
Then we continue the greedy algorithm to choose
these edges for the remainder of $A_0 \cup B_0$.
This is possible as the number of choices at each step 
is at least $|A'_1|/2$, of which at most $3|A_0|$ are forbidden
due to using a colour or a vertex used at a previous step.

Finally, consider $G_3=G[A_3,B_3]$ obtained from $G'_1$ 
by deleting all vertices covered by $M_0$
and all edges that share a colour with $M_0$.
Then $|A_3|=|B_3|=|A'_1|-|A_0|>(1-10^{-3}d)|A_1|$
and $G_3$ has minimum degree at least
$10^{-3}d|A_1| - |A_0| - |M_0| \ge 10^{-3} d|A_1|/2$.

We claim that $G_3$ has a perfect matching.
To see this we check Hall's condition.
Suppose for a contradiction there is $S \sub A_3$ with $|N(S)|<|S|$. 
By the minimum degree we have $|S| \ge 10^{-3}d|A_1|/2$.
Now we cannot have $|S| \le |A_1|/2$, as then
$|N_{G_3}(S)| \ge \min\{2|S|,2|A_1|/3\} - 2|A_0| > |S|$.
However, letting $T = B_3 \sm N(S)$ we have $N(T) \sub A_3 \sm S$,
so $|N(T)|=|A_3|-|S|<|B_3|-|N(S)|=|T|$. The same argument
as for $S$ gives $|T|>|A_1|/2$, contradiction.
Therefore $G_3$ has a perfect matching $M_3$.
Now $M_2 \cup M_0 \cup M_3$ is a rainbow perfect matching in $G$,
which completes the proof of Theorem \ref{main}.

\section{A robustly matchable pair} 

In this section we prove Lemma \ref{robmatch}.
Let $G$ be a bipartite graph with parts $A$ and $B$.
We say $G$ is $(\eps,\dD)$-dense if 
for any $A' \sub A$ and $B' \sub B$
with $|A'| \ge \eps |A|$ and $|B'| \ge \eps |B|$
we have $e_G(A',B') \ge \dD |A'||B'|$.
We start by applying the following
result of Peng, R\"odl and Ruci\'nski \cite[Theorem 1.3]{PRR}
with $\eps=1/10$, $c=0.24$ and $c'=1/50$.%
\footnote{
This result follows from their proof; 
they state the case $c=1/8$, $c'=1/2$ 
and use $\log_2(1+\eps/8) \ge \eps/6$.}

\begin{lemma} 
Suppose $c,c' \in (0,1)$ with $4c+c'\leq 1$.
Then there are $A' \sub A$ and $B' \sub B$
with $|A'|=|B'| \ge d^{2/\log_2(1+c\eps)} n/2$
so that $G'=G[A',B']$ is $(\eps,c'd)$-dense.
\end{lemma}

Let $G_1=G[A_1,B_1]$ be obtained from $G'$ 
by the following algorithm. 
Initially, $A_1=A'$ and $B_1=B'$.
At any step of the algorithm, we update $G_1$ 
by deleting a vertex or set of vertices 
of one of the following types 
(choosing arbitrarily if there is a choice).
\begin{enumerate}
\item $v \in A_1$ or $v \in B_1$ with $d_{G_1}(v) \le 10^{-3} d|A'|$,
\item $S \sub A_1$ or $S \sub B_1$ with
$|S|<\eps|A'|$ and $|N_{G_1}(S)| \le 2|S|$.
\end{enumerate}
Whenever we delete some vertices from $A_1$ or $B_1$
we delete an arbitrary set of the same size from the other,
so that we always maintain $|A_1|=|B_1|$.
We stop if no deletion is possible or if we have 
deleted at least $2\eps|A'|$ vertices from each side.

We claim that the latter option is impossible. 
Indeed, then without loss of generality we deleted
$\eps |A'|$ vertices from $A'$ of type (i) or (ii) as above 
(at least half of the deleted vertices are deleted
for a reason other than maintaining equal part sizes).
Let $D^A = D^A_i \cup D^A_{ii}$ be the deleted vertices
in $A'$ according to deletions of type (i) or (ii).
Note that $|D^A|< 3\eps|A'|$, and
$|N_{G_1}(D^A_{ii})| \le 2|D^A_{ii}| < 6\eps|A'|$.
Let $B_0 = B_1 \sm N_{G_1}(D^A_{ii})$, so $|B_0| > (1-9\eps)|A'| = \eps|A'|$.
Now $e_{G'}(D^A,B_0) \le |D^A_i| \cdot 10^{-3}d|A'| < \tfrac{d}{50} |D^A||B_0|$
contradicts $(\eps,d/50)$-density of $G'$, which proves the claim.

Thus the algorithm stops with $|A_1|=|B_1|>(1-3\eps)|A'| \ge d^{60} n/3$
(using $2/\log_2(1.024) < 60$),
minimum degree at least $10^{-3}d|A'|$ and $|N_{G_1}(S)| \ge 2|S|$
for any $S \sub A_1$ or $S \sub B_1$ with $|S|<\eps |A'|$.
Furthermore, for any $S \sub A_1$ or $S \sub B_1$
with $|S| \ge \eps |A'|$, by $(\eps,d/50)$-density of $G'$
we have $|N_{G_1}(S)| \ge |B_1|-\eps|A'| \ge 2|A_1|/3$.
This proves Lemma \ref{robmatch}.

\section{Augmentation algorithm} \label{sec:alg}

In this section we prove Lemma \ref{alg},
which will complete the proof of Theorem \ref{main}.
Suppose it is not true, say $|A'_0| > p|A_1|/4$.
We will iteratively construct 
$R = R^A \cup R^B \sub M_2$,
where we think of $R^A$ and $R^B$ as
`reachable' from $A_0$ and $B_0$.
At some point $R^A$ and $R^B$ will intersect,
which will contradict $M_2$ being 
a maximum size rainbow matching in $G^*$.
Let $\tT := n^{-0.66}$, and note that
$\tT |A_1| > n^{0.04}/3$.

\begin{alg}
Let $R^A = R^B = \es$ and let $C$ 
be the set of colours not used by $M_2$.
At step $i \ge 1$:
\begin{enumerate}
\item if $R^A \cap R^B \ne \es$ stop, otherwise
let $R^A_i$ be the set of all $uv \in M_2$
where $v \in B_2 \sm V(R^B)$ such that at least $\tT |A_1|$ edges 
in $G^*$ from $v$ to $A_0$ use a colour in $C$,
let $C^A_i$ be the set of colours used by $R^A_i$,
update $R^A$ by adding $R^A_i$ and $C$ by adding $C^A_i$,
\item if $R^A \cap R^B \ne \es$ stop, otherwise
let $R^B_i$ be the set of all $uv \in M_2$
where $u \in A_2 \sm V(R^A)$ such that at least $\tT |A_1|$ edges 
in $G^*$ from $u$ to $B_0$ use a colour in $C$,
let $C^B_i$ be the set of colours used by $R^B_i$,
update $R^B$ by adding $R^B_i$ and $C$ by adding $C^B_i$.
\end{enumerate}
\end{alg}

\begin{claim}
$|R^A_1| \ge |A_1|/4$.
\end{claim}

To see this, we consider the number $X$ of edges in $G^*$ with colour in $C$
between $A'_0$ and $B_2$. We have $X \le |R^A_1| |A'_0| + |B_2| \tT |A_1|$
by definition of $R^A_1$. Also, by definition of $A'_0$, every vertex in $A'_0$ 
has at least $(1-2p)|B_2| - (|M_2| - |B'_1|/2) \ge |A_1|/3$ edges in $G^*$ 
to $B_2$ with colour in $C$, so $X \ge |A'_0| \cdot |A_1|/3$.
As $|A'_0| \ge p|A_1|/4$ and $p^{-1} \tT n = n^{0.66} \ll pn \ll |A_1|$,
we deduce $|R^A_1| \ge |A_1|/3 - |A'_0|^{-1} |B_2| \tT |A_1|
\ge |A_1|/3 - 4p^{-1} \tT n \ge |A_1|/4$, as claimed.

\begin{claim}
For $i \ge 1$, we have $|R^B_i| \ge |R^A| - 3pn$
and $|R^A_{i+1}| \ge |R^B| - 3pn$.
\end{claim}

To see this, we first note that as $R^A \cap R^B = \es$, any vertex 
in $B_0$ has at least $(1-2p)|A_2| - |R^B| - |M_2 \sm (R^A \cup R^B)| 
= |R^A| - 2p|A_2|$ edges in $G^*$ to $A_2 \sm R^B$ with colour in $C$.
Double-counting such edges as in the previous claim gives
$|B_0|(|R^A| - 2p|A_2|) \le |R^B_i| |B_0| + |B_2| \tT |A_1|$,
so $|R^B_i| \ge |R^A| - 2p|A_2| - |B_0|^{-1} |B_2| \tT |A_1| \ge |R^A| - 3pn$.
The proof of the second inequality is similar, so the claim holds.

\begin{claim}
The algorithm terminates at some step $i = i^+ < \log n$.
\end{claim}

To see this, we show inductively that if $R^A \cap R^B = \es$ at step $i$
then $|R^B_i| \ge f(i) |A_1|/3$ and $|R^A_i| \ge (f(i)-2^{-3}) |A_1|/3$
where $f(i) = 2^{i-4} + 2^{-1}$. First note that $f(1) = 5/8$
and for $i \ge 2$ we have $\sum_{j=1}^{i-1} (f(j)-2^{-3}) 
= 2^{-4}(2^i-1)+(i-1)(2^{-1}-2^{-3}) \ge f(i) - 3/16$.
At step $1$ we have $|R^A_1| \ge |A_1|/4 > (f(1)-2^{-3})|A_1|/3$ 
and $|R^B_1| \ge |R^A| - 3pn > 0.21|A_1| > f(1)|A_1|/3$.
Supposing the statement at step $i-1 \ge 1$, 
we have $|R^A_i| \ge (\sum_{j=1}^{i-1} |R^B_j|)-3pn
\ge (\sum_{j=1}^{i-1} f(j) ) |A_1|/3 - 3pn
\ge (f(i) - 1/16 ) |A_1|/3 - 3pn \ge (f(i)-2^{-3}) |A_1|/3$
and $|R^B_i| \ge (\sum_{j=1}^i |R^A_j|)-3pn
\ge (\sum_{j=1}^i f(j) - 2^{-3}) |A_1|/3 - 3pn
\ge (f(i+1)-3/16) |A_1|/3 - 3pn \ge f(i) |A_1|/3$.
Thus the required bounds hold by induction.
While $R^A \cap R^B = \es$ we deduce
$(2f(i)-2^{-3}) |A_1|/3 < |M_2| < n$,
so $i^+ < \log n$, as claimed.

\medskip

The algorithm terminates by finding some edge 
$ab \in R^A \cap R^B$ where $a \in A_2$ and $b \in B_2$.
We will obtain a contradiction by modifying $M_2$ 
to obtain a larger rainbow matching in $G^*$.
Given two colours $c$ and $c'$ in $C$, 
we say that $c$ is earlier than $c'$ if 
$c$ was added to $C$ before $c'$.
We start by applying the definition of $R^A$ and $R^B$
to find edges $a_0 b$ and $a b_0$ of $G^*$ with 
$a_0 \in A_0$ and $b_0 \in B_0$ where the colours
of $a_0 b$ and $a b_0$ are in $C$ and earlier than that of $ab$.
We modify $M_2$ to obtain $M'_2$
by deleting $ab$ and adding $a_0b$ and $ab_0$.
Thus we obtain a larger matching, but $M'_2$ may not be rainbow,
due to repeating the colours of $a_0b$ and $ab_0$.
While the current matching $M'_2$ is not rainbow,
we apply the following `trace back' algorithm 
(similar to that of \cite{GS}).

\begin{alg} \label{alg}
At step $i \ge 1$ we have at most two `active' edges,
which are edges of $M'_2$ having some colour in $C$ 
shared with some edge that is still present from $M_2$.
At step $1$ these are $a_0b$ and $ab_0$.
If there is an active edge at step $i$, we choose one arbitarily, 
call it $a_i b_i$, and let $a'_i b'_i$ be the edge of $M_2$ of the same colour $c \in C$.
By construction of $C$, one of $a'_i$ or $b'_i$, say $a'_i$, has at least $\tT |A_1|$
edges to $B_0$ or $A_0$ using an earlier colour than $c$ in $C$.
We modify $M'_2$ by deleting $a'_ib'_i$ and adding some such edge $a'_i b^i_0$
where $b^i_0 \in B_0$ is distinct from all previous choices.
We say that $a_i b_i$ is no longer active.
We make $a'_i b^i_0$ active if its colour is shared with some edge 
that is still present from $M_2$.
\end{alg}

Algorithm \ref{alg} is illustrated in Figure \ref{fig:pf}:
the thick black edge represents the edge $ab \in R^A \cap R^B$,
at step 1 the green and blue diagonals are active,
at step 2 the blue diagonal is active,
at step 3 the red diagonal is active,
at step 4 the pink diagonal is active,
at step 5 there are no active edges
so the algorithm terminates.
To see that the algorithm succeeds,
note that there are at most $4\log n$ steps
of replacing an active edge by another,
and each choice has at least $\tT |A_1| > n^{0.04}/3 > 4\log n$ options.
Thus we obtain a rainbow matching $M'_2$ in $G^*$ with $|M'_2|>|M_2|$.
This contradiction proves Lemma \ref{alg}.

\medskip

\nib{Postscript.}
The Akbari-Alipour conjecture was proved independently and simultaneously
by Montgomery, Pokrovskiy and Sudakov \cite{MPS}.
Our proof is much simpler than theirs, and gives a better bound on the
number of symbols required, whereas their proof applies in a 
much more general setting, and so has several further applications.
Results similar to those in \cite{MPS} 
(but not including the Akbari-Alipour conjecture)
were independently and simultaneously
obtained by Kim, K\"uhn, Kupavskii and Osthus \cite{KKKO}.

\end{document}